\documentclass[12pt,reqno]{amsart}
\usepackage{amsfonts}
\usepackage{amssymb,amsmath}
\usepackage{amscd}
\usepackage{hyperref}
\usepackage{graphicx}
\usepackage{epsfig}

\newtheorem{theorem}{Theorem}[section]
\newtheorem{lemma}{Lemma}[section]

\newtheorem{definition}{Definition}[section]
\newtheorem{example}{Example}[section]

\numberwithin{equation}{section}
\begin{document}
\title[Sharp inequalities for anti-invariant Riemannian Submersions]{Sharp
inequalities for anti-invariant Riemannian Submersions from Sasakian Space
forms}
\author[H. AYT\.{I}MUR]{H\"{U}LYA\ AYT\.{I}MUR}
\address{Bal\i kesir University\\
Department of Mathematics\\
10145 Bal\i kesir, TURKEY}
\email{hulya.aytimur@balikesir.edu.tr}
\author[C. \"{O}ZG\"{U}R ]{C\.{I}HAN \"{O}ZG\"{U}R}
\address{Bal\i kesir University\\
Department of Mathematics\\
10145 Bal\i kesir, TURKEY}
\email{cozgur@balikesir.edu.tr}

\begin{abstract}
We obtain sharp inequalities involving the Ricci curvature and the scalar
curvature for anti-invariant Riemannian submersions from Sasakian space
forms onto Riemannian manifolds.
\end{abstract}

\maketitle

\section{\textbf{Introduction}\label{sec:intro}}

To find relationship between the extrinsic and intrinsic invariants of a
submanifold have been very popular problems in the recent twenty five years.
The first study in this direction was started by B.-Y. Chen in 1993. He
established some inequalities between the main extrinsic (the squared mean
curvature) and main intrinsic invariants (the scalar curvature and the Ricci
curvature) of a submanifold in a real space form \cite{Chen-93}. In 1999,
Chen also established a relation between the Ricci curvature and the squared
mean curvature for a submanifold \cite{Chen-99}. After that, many papers
have been published by various authors in different ambient spaces. In 2011,
Chen published a book which consists of the all studies doing in these
directions \cite{Chen-2011}. The topic is still very popular and there are
many new papers related to the inequalities which are introduced by Chen.
For example see \cite{AyOz}, \cite{Aydin-2015}, \cite{Chen-99}, \cite%
{Adela-06}, \cite{MM-18}, \cite{Mihai-Ozg-TJM}, \cite{Mihai-Ozgur-Rocky},
\cite{ozgur-2011} and \cite{Sahin-2016}.

Let $\left( M,g\right) $ and $\left( B,g^{^{\prime }}\right) $ be $m$ and $b$%
-dimensional Riemannian manifolds, respectively. A \textit{Riemannian
submersion }$\pi :M\rightarrow B$ is a mapping of $M$ onto $B$ such that $%
\pi $ has a maximal rank and the differential $\pi _{\ast }$ preserves the
lengths of the horizontal vectors \cite{O'Neill-66}. In \cite{Chen-05} and
\cite{Chen-05b}, Chen proved a simple optimal relationship between
Riemannian submersions and minimal immersions \cite{Chen-05}. In \cite{Pab},
Alegre, Chen and Munteanu established a sharp relationship between the $%
\delta $-invariants and Riemannian submersions with totally geodesic fibers.
In \cite{Gul-17}, G\"{u}lbahar, Meri\c{c} and K\i l\i \c{c} obtained sharp
inequalities involving the Ricci curvature for Riemannian submersions. In
\cite{Sahin-2010}, \c{S}ahin introduced anti-invariant Riemannian
submersions from almost Hermitian manifolds onto Riemannian manifolds.

Motivated by the above studies, in the present study, we consider
anti-invariant Riemannian submersions from Sasakian manifolds onto
Riemannian manifolds. We obtain sharp inequalities involving the Ricci
curvature and the scalar curvature.

The paper is organized as follows. In Section 2, we give brief introduction
about Sasakian manifolds and submersions. We give some lemmas which will be
used in Section \ref{sec:2} and Section \ref{sec:3}. In Section \ref{sec:2},
we obtain some inequalities involving the Ricci curvature and the scalar
curvature on the vertical and horizontal distributions for anti-invariant
Riemannian submersions from Sasakian space forms. The equality cases are
also discussed. In Section \ref{sec:3}, we prove Chen-Ricci inequalities on
the vertical and horizontal distributions for anti-invariant Riemannian
submersions from Sasakian space forms. We find relationships between the
intrinsic and extrinsic invariants using fundamental tensors. The equality
cases are also considered.

\section{\textbf{Preliminaries} \label{sec:1}}

Let $\pi :M\rightarrow B$ be a Riemannian submersion. We put $\dim M=2m+1$
and $\dim B=b.$ For $x\in B$, Riemannian submanifold $\pi ^{-1}\left(
x\right) $ with the induced metric $\overline{g}$ is called a \emph{fiber}
and denoted by $\overline{M}.$ We notice that the dimension of each fiber is
always $\left( 2m+1-b\right) =r$ and dimension of the horizontal
distribution is $n=\left( 2m+1-r\right) $. In the tangent bundle $TM$ of $M$%
, the vertical and horizontal distributions are denoted by $\mathcal{V}%
\left( M\right) $ and $\mathcal{H}\left( M\right) $, respectively. We call a
vector field $X$ on $M$ \emph{projectable} if there exists a vector field $%
X_{\ast }$ on $B$ such that $\pi _{\ast }\left( X_{p}\right) =X_{\ast \pi
\left( p\right) }$ for each $p\in M.$ In this case, we call that $X$ and $%
X_{\ast }$ are $\pi $-related. A vector field $X$ on $M$ is called \emph{%
basic} if it is projectable and horizontal (\cite{O'Neill-66} and \cite%
{O'Neill 83}).

The tensor fields $T$ and $A$ of type $\left( 1,2\right) $ are defined by
\begin{equation*}
T_{E}F=h\nabla _{\upsilon E}\upsilon F+\upsilon \nabla _{\upsilon E}hF,
\end{equation*}%
\begin{equation*}
A_{E}F=h\nabla _{hE}\upsilon F+\upsilon \nabla _{hE}hF.
\end{equation*}

Denote by $R,$ $R^{^{\prime }},$ $\widehat{R}$ and $R^{\ast }$ the
Riemannian curvature tensor of Riemannian manifolds $M,$ $B,$ the vertical
distribution $\mathcal{V}$ and the horizontal distribution $\mathcal{H}$,
respectively. Then the Gauss-Codazzi type equations are given by
\begin{equation}
R\left( U,V,F,W\right) =\widehat{R}\left( U,V,F,W\right) +g\left(
T_{U}W,T_{V}F\right) -g\left( T_{V}W,T_{U}F\right) ,  \label{gauss1}
\end{equation}%
\begin{equation*}
R\left( X,Y,Z,H\right) =R^{\ast }\left( X,Y,Z,H\right) -2g\left(
A_{X}Y,A_{Z}H\right)
\end{equation*}%
\begin{equation}
+g\left( A_{Y}Z,A_{X}H\right) -\left( A_{X}Z,A_{Y}H\right) ,  \label{gauss2}
\end{equation}%
\begin{equation*}
R\left( X,V,Y,W\right) =g\left( \left( \nabla _{X}T\right) \left( V,W\right)
,Y\right) +g\left( \left( \nabla _{V}A\right) \left( X,Y\right) ,W\right)
\end{equation*}%
\begin{equation}
-g\left( T_{V}X,T_{W}Y\right) +g\left( A_{Y}W,A_{X}V\right) ,  \label{gauss3}
\end{equation}%
where
\begin{equation*}
\pi _{\ast }\left( R^{\ast }\left( X,Y\right) Z\right) =R^{^{\prime }}\left(
\pi _{\ast }X,\pi _{\ast }Y\right) \pi _{\ast }Z
\end{equation*}%
for any $X,Y,Z,H\in \chi ^{h}\left( M\right) $ and $U,V,F,W\in \chi
^{v}\left( M\right) $ \cite{O'Neill-66}.

Moreover, the mean curvature vector field $H$ of any fibre of Riemannian
submersion $\pi $ is given by%
\begin{equation*}
H=rN,\text{ \ \ \ \ }N=\sum_{j=1}^{r}T_{U_{j}}U_{j}
\end{equation*}%
where $\left\{ U_{1},...,U_{r}\right\} $ is an orthonormal basis of the
vertical distribution $\mathcal{V}.$ Furthermore, $\pi $ \textit{has totally
geodesic fibers }if $T$ vanishes on $\chi ^{h}\left( M\right) $ and $\chi
^{v}\left( M\right) $.

Now we give the following lemmas:

\begin{lemma}
\label{lmma1}\cite{Falci-2004} Let $\left( M,g\right) $ and $\left(
B,g^{^{\prime }}\right) $ be Riemannian manifolds admitting a Riemannian
submersion $\pi :M\rightarrow B.$ For $E,F,G\in \chi \left( M\right) ,$ we
have%
\begin{equation*}
g\left( T_{E}F,G\right) =-g\left( F,T_{E}G\right) ,
\end{equation*}%
\begin{equation*}
g\left( A_{E}F,G\right) =-g\left( F,A_{E}G\right) .
\end{equation*}%
That is, $A_{E}$ and $T_{E}$ are anti-symmetric with respect to $g.$
\end{lemma}

\begin{lemma}
\label{lemm3}\cite{Falci-2004} Let $\left( M,g\right) $ and $\left(
B,g^{^{\prime }}\right) $ be Riemannian manifolds admitting a Riemannian
submersion $\pi :M\rightarrow B.$

$(i)$ For $U,V\in \chi ^{v}\left( M\right) ,$%
\begin{equation*}
T_{U}V=T_{V}U,
\end{equation*}

$(ii)$ For $X,Y\in \chi ^{h}\left( M\right) ,$%
\begin{equation*}
A_{X}Y=-A_{Y}X.
\end{equation*}
\end{lemma}

Let $M$ be a ($2m+1$)-dimensional manifold and $\phi ,\xi ,\eta $ a tensor
field of type $\left( 1,1\right) $, a vector field, a $1$-form on $M,$
respectively. If $\phi ,\xi $ and $\eta $ satisfy the following conditions%
\begin{equation*}
\eta \left( \xi \right) =1,\text{ \ \ \ \ \ \ \ \ \ \ \ \ \ \ }\phi
^{2}X=-X+\eta \left( X\right) \xi
\end{equation*}%
for $X\in TM$, then $M$ is said to have an \textit{almost contact structure}
$\left( \phi ,\xi ,\eta \right) $ and $\left( M,\phi ,\xi ,\eta \right) $ is
called an\textit{\ almost contact manifold}. If
\begin{equation*}
\nabla _{X}\xi =-\phi X,\text{ \ \ \ \ \ }\left( \nabla _{X}\phi \right)
Y=g\left( X,Y\right) \xi -\eta \left( Y\right) X,
\end{equation*}%
then $\left( M,\nabla ,g,\phi ,\xi ,\eta \right) $ is called a \textit{%
Sasakian manifold}\emph{\ }\cite{Blair}, where $\nabla $ denotes the
Levi-Civita connection of $g$. $\phi $ is anti-symmetric with respect to $g$%
, that is, for $X,Y\in TM$
\begin{equation*}
g\left( \phi X,Y\right) +g\left( X,\phi Y\right) =0.
\end{equation*}%
A plane section $\pi $ in $TM$ is called a $\phi $-\textit{section} if it is
spanned by $X$ and $\phi X$, where $X$ is a unit tangent vector field
orthogonal to $\xi $. The sectional curvature of a $\phi $-section is called
a $\phi $-\textit{sectional curvature}. A Sasakian manifold with constant $%
\phi $-sectional curvature $c$ is said to be a \textit{Sasakian space form }%
\cite{Blair} and is denoted by $M(c)$. The curvature tensor $R$ of $M(c)$ is
expressed by

\begin{equation*}
R(X,Y)Z=\frac{c+3}{4}[g(Y,Z)X-g(X,Z)Y]+\frac{c-1}{4}[\eta (X)\eta (Z)Y
\end{equation*}%
\begin{equation*}
-\eta (Y)\eta (Z)X+g(X,Z)\eta (Y)\xi -g(Y,Z)\eta (X)\xi +g(\phi Y,Z)\phi X
\end{equation*}%
\begin{equation}
-g(\phi X,Z)\phi Y-2g(\phi X,Y)\phi Z].  \label{curv1}
\end{equation}

\begin{definition}
\label{dfn}\cite{Kupeli-Murathan-2013} Let $\left( M,\nabla ,g,\phi ,\xi
,\eta \right) $ be a Sasakian manifold and $\left( B,g^{^{\prime }}\right) $
a Riemannian manifold. A Riemannian submersion $\pi :M\rightarrow B$ is
called anti-invariant if \ $\mathcal{V}\left( M\right) $ is anti-invariant
with respect to $\phi $, i.e. $\phi \left( \mathcal{V}\left( M\right)
\right) \subseteq \mathcal{H}\left( M\right) $.
\end{definition}

Let $\pi :\left( M,\nabla ,g,\phi ,\xi ,\eta \right) \rightarrow \left(
B,g^{^{\prime }}\right) $ be an anti-invariant Riemannian submersion from a
Sasakian manifold $\left( M,\nabla ,g,\phi ,\xi ,\eta \right) $ to a
Riemannian manifold $\left( B,g^{^{\prime }}\right) $. From Definition \ref%
{dfn}, we have $\phi \left( \mathcal{V}\left( M\right) \right) \cap \mathcal{%
H}\left( M\right) \neq \left\{ 0\right\} .$ We denote the complementary
orthogonal distribution to $\phi \left( \mathcal{V}\left( M\right) \right) $
in $\mathcal{H}\left( M\right) $ by $\mu .$ Then we have
\begin{equation*}
\mathcal{H}\left( M\right) =\phi \left( \mathcal{V}\left( M\right) \right)
\oplus \mu .
\end{equation*}

Suppose that $\xi $ is vertical. It is easy to see that $\mu $ is an
invariant distribution of $\mathcal{H}\left( M\right) $ under the
endomorphism $\phi .$ Thus for $X\in \chi ^{h}\left( M\right) ,$ we write%
\begin{equation*}
\phi X=BX+CX,
\end{equation*}%
where $BX\in \chi ^{v}\left( M\right) $ and $CX\in \chi \left( \mu \right) $
\cite{Kupeli-Murathan-2013}.

Suppose that $\xi $ is horizontal. It is easy to see that $\mu =\phi \mu
\oplus \left\{ \xi \right\} .$ Thus for $X\in \chi ^{h}\left( M\right) ,$ we
write%
\begin{equation*}
\phi X=BX+CX,
\end{equation*}%
where $BX\in \chi ^{v}\left( M\right) $ and $CX\in \chi \left( \mu \right) $
\cite{Kupeli-Murathan-2013}.

\begin{lemma}
\label{lemm2} \cite{Kupeli-Murathan-2013} Let $\pi :M\rightarrow B$ be an
anti-invariant Riemannian submersion from a Sasakian manifold $\left(
M,\nabla ,g,\phi ,\xi ,\eta \right) $ to a Riemannian manifold $\left(
B,g^{^{\prime }}\right) .$

$(i)$ If $\xi $ is vertical, then $C^{2}X=-X-\phi BX$

$(ii)$ If $\xi $ is horizontal, then $C^{2}X=-X+\eta \left( X\right) \xi
-\phi BX.$
\end{lemma}

\begin{example}
\label{exmp1}\cite{Blair} Let us take $M=%
\mathbb{R}
^{2m+1}$ with the standard coordinate functions $\left(
x_{1},...,x_{m},y_{1},...,y_{m},z\right) ,$ the contact structure $\eta =%
\frac{1}{2}(dz-\sum_{i=1}^{m}y_{i}dx_{i}),$ the characteristic vector field $%
\xi =2\frac{\partial }{\partial z}$ and the tensor field $\varphi $ given by%
\begin{equation*}
\varphi =%
\begin{bmatrix}
0 & \delta _{ij} & 0 \\
-\delta _{ij} & 0 & 0 \\
0 & y_{j} & 0%
\end{bmatrix}%
.
\end{equation*}%
The Riemannian metric is $g=\eta \otimes \eta +\frac{1}{4}\overset{m}{%
\underset{i=1}{\sum }}\left( (dx_{i})^{2}+(dy_{i})^{2}\right) .$ Then $%
\left( M^{2m+1},\varphi ,\xi ,\eta ,g\right) $ is a Sasakian space form with
constant $\varphi -$sectional curvature $c=-3$ and it is denoted by $%
\mathbb{R}
^{2m+1}(-3).$ The vector fields%
\begin{equation*}
E_{i}=2\frac{\partial }{\partial y_{i}},\text{ }E_{i+m}=\varphi X_{i}=2(%
\frac{\partial }{\partial x_{i}}+y_{i}\frac{\partial }{\partial z}),\text{ }%
1\leq i\leq m,\text{ }\xi =2\frac{\partial }{\partial z},
\end{equation*}%
form a $g$-orthonormal basis for the contact metric structure.
\end{example}

\begin{example}
\cite{Kupeli-Murathan-2013}We consider $M=%
\mathbb{R}
^{5}(-3)$ with the structure given in Example \ref{exmp1}. The Riemannian
metric $g_{\mathbb{R}^{2}}$ is given by%
\begin{equation*}
g_{\mathbb{R}^{2}}=\frac{1}{8}\left[
\begin{array}{cc}
1 & 0 \\
0 & 1%
\end{array}%
\right]
\end{equation*}%
on $\mathbb{R}^{2}$. Let $\pi :%
\mathbb{R}
^{5}(-3)\rightarrow \mathbb{R}^{2}$ be a map defined by
\begin{equation*}
\pi \left( x_{1},x_{2},y_{1},y_{2},z\right) =\left(
x_{1}+y_{1},x_{2}+y_{2}\right) .
\end{equation*}%
Then
\begin{equation*}
\mathcal{V}\left( M\right) =sp\left\{
V_{1}=E_{1}-E_{3},V_{2}=E_{2}-E_{4},V_{3}=E_{5}=\xi \right\}
\end{equation*}%
and%
\begin{equation*}
\mathcal{H}\left( M\right) =sp\left\{
H_{1}=E_{1}+E_{3},H_{2}=E_{2}+E_{4}\right\} .
\end{equation*}%
So $\pi $ is a Riemannian submersion. Moreover, $\phi V_{1}=H_{1},\phi
V_{2}=H_{2},\phi V_{3}=0$ imply that $\phi (\mathcal{V}\left( M\right) )=%
\mathcal{H}\left( M\right) .$ Hence $\pi $ is an anti-invariant Riemannian
submersion such that $\xi $ is vertical.
\end{example}

\begin{example}
\cite{Kupeli-Murathan-2013}We consider $M=%
\mathbb{R}
^{5}(-3)$ with the structure given in Example \ref{exmp1}. Let $N=%
\mathbb{R}
^{3}-\left\{ \left( y_{1},y_{2},z\right) \in
\mathbb{R}
^{3}\mid y_{1}^{2}+y_{2}^{2}\leq 2\right\} .$ The Riemannian metric tensor $%
g_{N}$ is given by
\begin{equation*}
g_{N}=\frac{1}{4}\left[
\begin{array}{ccc}
\frac{1}{2} & \frac{y_{1}y_{2}}{2} & -\frac{y_{1}}{2} \\
\frac{y_{1}y_{2}}{2} & \frac{1}{2} & -\frac{y_{2}}{2} \\
-\frac{y_{1}}{2} & -\frac{y_{2}}{2} & 1%
\end{array}%
\right]
\end{equation*}%
on $N.$ Let $\pi :%
\mathbb{R}
^{5}(-3)\rightarrow N$ be a map defined by
\begin{equation*}
\pi \left( x_{1},x_{2},y_{1},y_{2},z\right) =\left( x_{1}+y_{1},x_{2}+y_{2},%
\frac{y_{1}^{2}}{2}+\frac{y_{2}^{2}}{2}+z\right) .
\end{equation*}%
Then
\begin{equation*}
\mathcal{V}\left( M\right) =sp\left\{
V_{1}=E_{1}-E_{3},V_{2}=E_{2}-E_{4}\right\}
\end{equation*}%
and%
\begin{equation*}
\mathcal{H}\left( M\right) =sp\left\{
H_{1}=E_{1}+E_{3},H_{2}=E_{2}+E_{4},H_{3}=E_{5}=\xi \right\} .
\end{equation*}%
So $\pi $ is a Riemannian submersion. Moreover, $\phi V_{1}=H_{1},\phi
V_{2}=H_{2}$ imply that $\phi (\mathcal{V}\left( M\right) )\subset \mathcal{H%
}\left( M\right) =\phi (\mathcal{V}\left( M\right) )\oplus \left\{ \xi
\right\} .$ Hence $\pi $ is an anti-invariant Riemannian submersion such
that $\xi $ is horizontal.
\end{example}

\section{\textbf{Inequalities for anti-invariant Riemannian submersions }
\label{sec:2}}

In the present section, we aim to obtain some inequalities involving the
Ricci curvature and the scalar curvature on the vertical and horizontal
distributions for anti-invariant Riemannian submersions from Sasakian space
forms. We shall also consider the equality cases of these inequalities.

Let $\left( M(c),g\right) ,$ $\left( B,g^{^{\prime }}\right) $ be a Sasakian
space form and a Riemannian manifold, respectively and $\pi :M(c)\rightarrow
B$ an anti-invariant Riemannian submersion. Furthermore, let $\left\{
U_{1},...,U_{r},X_{1},...,X_{n}\right\} $ be an orthonormal basis of $TM(c)$
such that $\mathcal{V}$ $=span\left\{ U_{1},...,U_{r}\right\} ,$ $\mathcal{H=%
}span\left\{ X_{1},...,X_{n}\right\} $. Then using (\ref{curv1}) and (\ref%
{gauss1}), we have\newpage
\begin{equation*}
\widehat{R}\left( U,V,F,W\right) =\frac{c+3}{4}\left\{ g\left( V,F\right)
g\left( U,W\right) -g\left( U,F\right) g\left( V,W\right) \right\}
\end{equation*}%
\begin{equation*}
+\frac{c-1}{4}\left\{ \eta \left( U\right) \eta \left( F\right) g\left(
V,W\right) -\eta \left( V\right) \eta \left( F\right) g\left( U,W\right)
\right.
\end{equation*}%
\begin{equation*}
+\eta \left( V\right) \eta \left( W\right) g\left( U,F\right) -\eta \left(
U\right) \eta \left( W\right) g\left( V,F\right) +g\left( \phi V,F\right)
g\left( \phi U,W\right)
\end{equation*}%
\begin{equation*}
\left. -g\left( \phi V,W\right) g\left( \phi U,F\right) -2g\left( W,\phi
F\right) g\left( \phi U,V\right) \right\}
\end{equation*}%
\begin{equation}
-g\left( T_{U}W,T_{V}F\right) +g\left( T_{V}W,T_{U}F\right) .  \label{R1}
\end{equation}%
Similarly, from (\ref{curv1}) and (\ref{gauss2}), we get
\begin{equation*}
R^{\ast }\left( X,Y,Z,H\right) =\frac{c+3}{4}\left\{ g\left( Y,Z\right)
g\left( X,H\right) -g\left( X,Z\right) g\left( Y,H\right) \right\}
\end{equation*}%
\begin{equation*}
+\frac{c-1}{4}\left\{ \eta \left( X\right) \eta \left( Z\right) g\left(
Y,H\right) -\eta \left( Y\right) \eta \left( Z\right) g\left( X,H\right)
\right.
\end{equation*}%
\begin{equation*}
+\eta \left( Y\right) \eta \left( H\right) g\left( X,Z\right) -\eta \left(
X\right) \eta \left( H\right) g\left( Y,Z\right) +g\left( \phi Y,Z\right)
g\left( \phi X,H\right)
\end{equation*}%
\begin{equation*}
\left. -g\left( \phi Y,H\right) g\left( \phi X,Z\right) -2g\left( H,\phi
Z\right) g\left( \phi X,Y\right) \right\}
\end{equation*}%
\begin{equation}
+2g\left( A_{X}Y,A_{Z}H\right) -g\left( A_{Y}Z,A_{X}H\right) +\left(
A_{X}Z,A_{Y}H\right) .  \label{R2}
\end{equation}
\textbf{Case I:} Assume that $\xi $\ is vertical.

For the vertical distribution, in view of (\ref{R1}), since $\pi $ is
anti-invariant and $\xi $ is vertical, we find%
\begin{equation*}
\widehat{Ric}\left( U\right) =\frac{c+3}{4}\left( r-1\right) g\left(
U,U\right) +\frac{c-1}{4}\left\{ \left( 2-r\right) \eta \left( U\right)
^{2}-g\left( U,U\right) \right\}
\end{equation*}%
\begin{equation*}
-rg\left( T_{U}U,H\right) +\sum_{j=1}^{r}g\left(
T_{U_{j}}U,T_{U}U_{j}\right) .
\end{equation*}%
Hence we obtain the following theorem:

\begin{theorem}
\label{thrm1} Let $\pi :M(c)\rightarrow B$ be an anti-invariant Riemannian
submersion from a Sasakian space form $\left( M(c),g\right) $ onto a
Riemannian manifold $\left( B,g^{^{\prime }}\right) $ such that $\xi $\ is
vertical. Then%
\begin{equation*}
\widehat{Ric}\left( U\right) \geq \frac{c+3}{4}\left( r-1\right) -\frac{c-1}{%
4}\left\{ \left( r-2\right) \eta \left( U\right) ^{2}+1\right\} -rg\left(
T_{U}U,H\right) .
\end{equation*}%
The equality case of \ the inequality holds for a unit vertical vector field
$U\in \chi ^{\mathcal{V}}\left( M(c)\right) $ if and only if each fiber is
totally geodesic.
\end{theorem}

Similarly in view of (\ref{R1}), using the symmetry of $T$, we have
\begin{equation*}
2\widehat{\tau }=\frac{c+3}{4}r\left( r-1\right) +\frac{c-1}{4}\left(
2-2r\right) -r^{2}\left\Vert H\right\Vert ^{2}+\sum_{i,j=1}^{r}g\left(
T_{U_{i}}U_{j},T_{U_{i}}U_{j}\right) ,
\end{equation*}%
where $\widehat{\tau }=\underset{1\leq i<j\leq r}{\sum }\widehat{R}\left(
U_{i},U_{j},U_{j},U_{i}\right) .$ Then we can write%
\begin{equation*}
2\widehat{\tau }\geq \frac{c+3}{4}r\left( r-1\right) -\frac{c-1}{2}%
(r-1)-r^{2}\left\Vert H\right\Vert ^{2}.
\end{equation*}%
The equality case of the inequality holds if and only if $T=0$, which means
that each fiber is totally geodesic. Thus we can state the following theorem:

\begin{theorem}
\label{thrm2} Let $\pi :M(c)\rightarrow B$ be an anti-invariant Riemannian
submersion from a Sasakian space form $\left( M(c),g\right) $ onto a
Riemannian manifold $\left( B,g^{^{\prime }}\right) $ such that $\xi $\ is
vertical. Then%
\begin{equation*}
2\widehat{\tau }\geq \frac{c+3}{4}r\left( r-1\right) -\frac{c-1}{2}%
(r-1)-r^{2}\left\Vert H\right\Vert ^{2}.
\end{equation*}%
The equality case of the inequality holds if and only if each fiber is
totally geodesic.
\end{theorem}

For the horizontal distribution, in view of (\ref{R2}), since $\pi $ is
anti-invariant and $\xi $ is vertical, using the anti-symmetry of $A,$ we
find

\begin{equation*}
2\tau ^{\ast }=\frac{c+3}{4}n\left( n-1\right)
\end{equation*}%
\begin{equation}
+\sum_{i,j=1}^{n}\left[ \frac{3\left( c-1\right) }{4}g\left(
CX_{i},X_{j}\right) g\left( CX_{i},X_{j}\right) -3g\left(
A_{X_{i}}X_{j},A_{X_{i}}X_{j}\right) \right] .  \label{Ric2}
\end{equation}%
By the use of Lemma \ref{lemm2}, we obtain%
\begin{equation*}
2\tau ^{\ast }=\frac{c+3}{4}n\left( n-1\right) +\frac{3}{4}\left( c-1\right)
\left( n+tr\left( \phi B\right) \right) -\sum_{i,j=1}^{n}3g\left(
A_{X_{i}}X_{j},A_{X_{i}}X_{j}\right) .
\end{equation*}

Then we can write
\begin{equation}
2\tau ^{\ast }\leq \frac{c+3}{4}n\left( n-1\right) +\frac{3}{4}\left(
c-1\right) \left( n+tr\left( \phi B\right) \right) ,  \label{T2}
\end{equation}%
where $\tau ^{\ast }=\underset{1\leq i<j\leq n}{\sum }R^{\ast }\left(
X_{i},X_{j},X_{j},X_{i}\right) .$ The equality case of (\ref{T2}) holds if
and only if $A=0$, which means that the horizontal distribution is
integrable. So we can state the following theorem:

\begin{theorem}
\label{thrm3} Let $\pi :M(c)\rightarrow B$ be an anti-invariant Riemannian
submersion from a Sasakian space form $\left( M(c),g\right) $ onto a
Riemannian manifold $\left( B,g^{^{\prime }}\right) $ such that $\xi $\ is
vertical. Then%
\begin{equation*}
2\tau ^{\ast }\leq \frac{c+3}{4}n\left( n-1\right) +\frac{3}{4}\left(
c-1\right) \left( n+tr\left( \phi B\right) \right) .
\end{equation*}%
The equality case of (\ref{T2}) holds if and only if $\ \mathcal{H}(M)$ is
integrable.
\end{theorem}

\textbf{Case II: }Assume that $\xi $\ is horizontal.

From (\ref{R1}), since $\pi $ is anti-invariant submersion, after some
computations, we have%
\begin{equation*}
2\widehat{\tau }=\frac{c+3}{4}r\left( r-1\right) -r^{2}\left\Vert
H\right\Vert ^{2}+\sum_{i,j=1}^{r}g\left(
T_{U_{i}}U_{j},T_{U_{i}}U_{j}\right) .
\end{equation*}%
Hence we can state the following theorem:

\begin{theorem}
\label{thrm21} Let $\pi :M(c)\rightarrow B$ be an anti-invariant Riemannian
submersion from a Sasakian space form $\left( M(c),g\right) $ onto a
Riemannian manifold $\left( B,g^{^{\prime }}\right) $ such that $\xi $\ is
horizontal. Then%
\begin{equation*}
2\widehat{\tau }\geq \frac{c+3}{4}r\left( r-1\right) -r^{2}\left\Vert
H\right\Vert ^{2}.
\end{equation*}

The equality case of the inequality holds if and only if each fiber is
totally geodesic.
\end{theorem}

For the horizontal distribution, from (\ref{R2}), since $\xi $ is horizontal
and $A$ is anti-symmetric, after some computations, we have

\begin{equation*}
2\tau ^{\ast }=\frac{c+3}{4}n\left( n-1\right) +\sum_{i,j=1}^{n}\left[ \frac{%
c-1}{4}\left\{ 2-2n+3g\left( CX_{i},X_{j}\right) g\left( CX_{i},X_{j}\right)
\right\} \right.
\end{equation*}%
\begin{equation*}
\left. -3g\left( A_{X_{i}}X_{j},A_{X_{i}}X_{j}\right) \right] .
\end{equation*}%
Then using Lemma \ref{lemm2}, we obtain%
\begin{equation*}
2\tau ^{\ast }=\frac{c+3}{4}n\left( n-1\right) +\frac{c-1}{4}\left( 3tr\phi
B+n-1\right)
\end{equation*}%
\begin{equation*}
-\sum_{i,j=1}^{n}3g\left( A_{X_{i}}X_{j},A_{X_{i}}X_{j}\right) ,
\end{equation*}%
where $\tau ^{\ast }=\underset{1\leq i<j\leq n}{\sum }R^{\ast }\left(
X_{i},X_{j},X_{j},X_{i}\right) .$

So we can state the following theorem:

\begin{theorem}
Let $\pi :M(c)\rightarrow B$ be an anti-invariant Riemannian submersion from
a Sasakian space form $\left( M(c),g\right) $ onto a Riemannian manifold $%
\left( B,g^{^{\prime }}\right) $ such that $\xi $\ is horizontal. Then
\begin{equation*}
2\tau ^{\ast }\leq \frac{c+3}{4}n\left( n-1\right) +\frac{\left( c-1\right)
}{4}\left( 3tr\left( \phi B\right) +n-1\right) .
\end{equation*}
The equality case of the inequality holds if and only if $\mathcal{H}(M)$ is
integrable.
\end{theorem}

\section{\textbf{Chen-Ricci inequalities for anti-invariant Riemannian
submersions }\label{sec:3}}

In the present section, we aim to obtain Chen-Ricci inequality on the
vertical and horizontal distributions for anti-invariant Riemannian
submersions from a Sasakian space forms onto a Riemannian manifold. The
equality cases will be also considered.

Let $\left( M(c),g\right) $ be a Sasakian space form and $\left(
B,g^{^{\prime }}\right) $ a Riemannian manifold. Assume that $\pi
:M(c)\rightarrow B$ is an anti-invariant Riemannian submersion and $\left\{
U_{1},...,U_{r},X_{1},...,X_{n}\right\} $ is an orthonormal basis of $TM(c)$
such that $\mathcal{V}$ $=span\{U_{1}$,..., $\ U_{r}\}$, $\mathcal{H=}%
span\left\{ X_{1},...,X_{n}\right\} .$ Now we denote $T_{ij}^{s}$ by
\begin{equation}
T_{ij}^{s}=g\left( T_{Ui}U_{j},X_{s}\right) ,  \label{T3}
\end{equation}%
where $1\leq i,j\leq r$ and $1\leq s\leq n$ (see \cite{Gul-17}).

Similarly, we denote $A_{ij}^{\alpha }$ by%
\begin{equation}
A_{ij}^{\alpha }=g\left( A_{Xi}X_{j},U_{\alpha }\right) .  \label{A1}
\end{equation}%
where $1\leq i,j\leq n$ and $1\leq \alpha \leq r.$ From \cite{Gul-17}, we use%
\begin{equation}
\delta \left( N\right) =\underset{i=1}{\overset{n}{\sum }}\underset{k=1}{%
\overset{r}{\sum }}g\left( \left( \nabla _{X_{i}}T\right)
_{U_{k}}U_{k},X_{i}\right) .  \label{deltaN}
\end{equation}

\textbf{Case I: }Assume that $\xi $\ is vertical.

Then from (\ref{R1}), we have
\begin{equation*}
2\widehat{\tau }=\frac{c+3}{4}r\left( r-1\right) -\frac{c-1}{2}\left(
r-1\right) -r^{2}\left\Vert H\right\Vert ^{2}+\sum_{i,j=1}^{r}g\left(
T_{U_{i}}U_{j},T_{U_{i}}U_{j}\right) .
\end{equation*}%
\ Using (\ref{T3}) in the last equality and the symmetry of $T$, we can
write
\begin{equation}
2\widehat{\tau }=\frac{c+3}{4}r\left( r-1\right) -\frac{c-1}{2}\left(
r-1\right) -r^{2}\left\Vert H\right\Vert
^{2}+\sum_{s=1}^{n}\sum_{i,j=1}^{r}\left( T_{ij}^{s}\right) ^{2}.  \label{T4}
\end{equation}%
For a local orthonormal frame $\left\{ X_{i},U_{j}\right\} _{1\leq i\leq
n,1\leq j\leq r}$ on $M(c)$, such that the horizontal and vertical
distributions are spanned by $\left\{ X_{i},\right\} _{1\leq i\leq n}$ and $%
\left\{ U_{j}\right\} _{1\leq j\leq r},$ respectively, we know from \cite%
{Gul-17} that%
\begin{equation*}
\sum_{s=1}^{n}\sum_{i,j=1}^{r}\left( T_{ij}^{s}\right) ^{2}=\frac{1}{2}%
r^{2}\left\Vert H\right\Vert ^{2}+\frac{1}{2}\sum_{s=1}^{n}\left[
T_{11}^{s}-T_{22}^{s}-...-T_{rr}^{s}\right] ^{2}
\end{equation*}%
\begin{equation}
+2\sum_{s=1}^{n}\sum_{j=2}^{r}\left( T_{1j}^{s}\right)
^{2}-2\sum_{s=1}^{n}\sum_{2\leq i<j\leq r}^{r}\left[ T_{ii}^{s}T_{jj}^{s}-%
\left( T_{ij}^{s}\right) ^{2}\right] .  \label{T1}
\end{equation}%
So using the above equality in (\ref{T4}), we get%
\begin{equation*}
2\widehat{\tau }=\frac{c+3}{4}r\left( r-1\right) -\frac{c-1}{2}\left(
r-1\right)
\end{equation*}%
\begin{equation*}
-\frac{1}{2}r^{2}\left\Vert H\right\Vert ^{2}+\frac{1}{2}\sum_{s=1}^{n}\left[
T_{11}^{s}-T_{22}^{s}-...-T_{rr}^{s}\right] ^{2}
\end{equation*}%
\begin{equation*}
+2\sum_{s=1}^{n}\sum_{j=2}^{r}\left( T_{1j}^{s}\right)
^{2}-2\sum_{s=1}^{n}\sum_{2\leq i<j\leq r}^{r}\left[ T_{ii}^{s}T_{jj}^{s}-%
\left( T_{ij}^{s}\right) ^{2}\right] .
\end{equation*}%
Then from the last equality, we have%
\begin{equation*}
2\widehat{\tau }\geq \frac{c+3}{4}r\left( r-1\right) -\frac{c-1}{2}\left(
r-1\right)
\end{equation*}%
\begin{equation}
-\frac{1}{2}r^{2}\left\Vert H\right\Vert ^{2}-2\sum_{s=1}^{n}\sum_{2\leq
i<j\leq r}^{r}\left[ T_{ii}^{s}T_{jj}^{s}-\left( T_{ij}^{s}\right) ^{2}%
\right] .  \label{T5}
\end{equation}%
Furthermore, from (\ref{gauss1}), taking $U=W=U_{i},$ $V=F=U_{j}$ and using (%
\ref{T3}) we can write%
\begin{equation*}
2\sum_{2\leq i<j\leq r}R\left( U_{i},U_{j},U_{j},U_{i}\right) =2\sum_{2\leq
i<j\leq r}\widehat{R}\left( U_{i},U_{j},U_{j},U_{i}\right)
\end{equation*}%
\begin{equation*}
+2\sum_{s=1}^{n}\sum_{2\leq i<j\leq r}\left[ T_{ii}^{s}T_{jj}^{s}-\left(
T_{ij}^{s}\right) ^{2}\right] .
\end{equation*}%
In view of the last equality, (\ref{T5}) can be written as%
\begin{equation*}
2\widehat{\tau }\geq \frac{c+3}{4}r\left( r-1\right) -\frac{c-1}{2}\left(
r-1\right) -\frac{1}{2}r^{2}\left\Vert H\right\Vert ^{2}
\end{equation*}%
\begin{equation}
+2\sum_{2\leq i<j\leq r}\widehat{R}\left( U_{i},U_{j},U_{j},U_{i}\right)
-2\sum_{2\leq i<j\leq r}R\left( U_{i},U_{j},U_{j},U_{i}\right) .  \label{z1}
\end{equation}%
Then using the equality
\begin{equation}
2\widehat{\tau }=2\underset{2\leq i<j\leq r}{\sum }\widehat{R}\left(
U_{i},U_{j},U_{j},U_{i}\right) +2\underset{j=1}{\overset{r}{\sum }}\widehat{R%
}\left( U_{1},U_{j},U_{j},U_{1}\right) ,  \label{T1a}
\end{equation}%
in view of (\ref{z1}), we have
\begin{equation*}
2\widehat{Ric}\left( U_{1}\right) \geq \frac{c+3}{4}r\left( r-1\right) -%
\frac{c-1}{2}\left( r-1\right)
\end{equation*}%
\begin{equation*}
-\frac{1}{2}r^{2}\left\Vert H\right\Vert ^{2}-2\sum_{2\leq i<j\leq r}R\left(
U_{i},U_{j},U_{j},U_{i}\right) .
\end{equation*}%
Since $M$ is a Sasakian space form, its curvature tensor $R$ satisfies the
equality (\ref{curv1}). So we obtain%
\begin{equation*}
\widehat{Ric}\left( U_{1}\right) \geq \frac{c+3}{4}\left( r-1\right) +\frac{%
c-1}{4}\left\{ (2-r)\eta \left( U_{1}\right) ^{2}-1\right\} -\frac{1}{4}%
r^{2}\left\Vert H\right\Vert ^{2}.
\end{equation*}%
Hence we state the following theorem:

\begin{theorem}
\label{thrm4} Let $\pi :M(c)\rightarrow B$ be an anti-invariant Riemannian
submersion from a Sasakian space form $\left( M(c),g\right) $ onto a
Riemannian manifold $\left( B,g^{^{\prime }}\right) $ such that $\xi $\ is
vertical. Then%
\begin{equation*}
\widehat{Ric}\left( U_{1}\right) \geq \frac{c+3}{4}\left( r-1\right) -\frac{%
c-1}{4}\left\{ (r-2)\eta \left( U_{1}\right) ^{2}+1\right\} -\frac{1}{4}%
r^{2}\left\Vert H\right\Vert ^{2}.
\end{equation*}%
The equality case of the inequality holds if and only if%
\begin{equation*}
T_{11}^{s}=T_{22}^{s}+...+T_{rr}^{s},
\end{equation*}%
\begin{equation*}
T_{1j}=0,\text{ \ }j=2,...,r.
\end{equation*}
\end{theorem}

On the other hand, using (\ref{A1}) and Lemma \ref{lemm2}, the equation (\ref%
{Ric2}) can be rewritten as%
\begin{equation*}
2\tau ^{\ast }=\frac{c+3}{4}n\left( n-1\right) +\frac{3}{4}(c-1)\left(
n+tr\left( \phi B\right) \right) -3\sum_{\alpha
=1}^{r}\sum_{i,j=1}^{n}\left( A_{ij}^{\alpha }\right) ^{2}.
\end{equation*}%
Since $A$ is anti-symmetric on $\chi ^{\mathcal{H}}\left( M(c)\right) ,$ the
above equality turns into%
\begin{equation*}
2\tau ^{\ast }=\frac{c+3}{4}n\left( n-1\right) +\frac{3}{4}(c-1)\left(
n+tr\left( \phi B\right) \right)
\end{equation*}%
\begin{equation}
-6\sum_{\alpha =1}^{r}\sum_{j=2}^{n}\left( A_{1j}^{\alpha }\right)
^{2}-6\sum_{\alpha =1}^{r}\sum_{2\leq i<j\leq n}\left( A_{ij}^{\alpha
}\right) ^{2}.  \label{A2}
\end{equation}%
Furthermore, from (\ref{gauss2}), taking $X=H=X_{i},$ $Y=Z=X_{j}$ and using (%
\ref{A1}), we have%
\begin{equation*}
2\sum_{2\leq i<j\leq n}R\left( X_{i},X_{j},X_{j},X_{i}\right) =2\sum_{2\leq
i<j\leq n}R^{\ast }\left( X_{i},X_{j},X_{j},X_{i}\right)
\end{equation*}%
\begin{equation}
+6\sum_{\alpha =1}^{r}\sum_{2\leq i<j\leq n}\left( A_{ij}^{\alpha }\right)
^{2}.  \label{T1aa}
\end{equation}%
If we consider the last equality in (\ref{A2}), then we get
\begin{equation*}
2\tau ^{\ast }=\frac{c+3}{4}n\left( n-1\right) +\frac{3}{4}(c-1)\left(
n+tr\left( \phi B\right) \right) -6\sum_{\alpha =1}^{r}\sum_{j=2}^{n}\left(
A_{1j}^{\alpha }\right) ^{2}
\end{equation*}%
\begin{equation*}
+2\sum_{2\leq i<j\leq n}R^{\ast }\left( X_{i},X_{j},X_{j},X_{i}\right)
-2\sum_{2\leq i<j\leq n}R\left( X_{i},X_{j},X_{j},X_{i}\right) .
\end{equation*}%
Since $M$ is a Sasakian space form, its curvature tensor $R$ satisfies the
equality (\ref{curv1}). Then we have
\begin{equation*}
2Ric^{\ast }\left( X_{1}\right) =\frac{c+3}{2}\left( n-1\right) +\frac{3}{4}%
(c-1)\left\Vert CX_{1}\right\Vert ^{2}
\end{equation*}%
\begin{equation*}
-6\sum_{\alpha =1}^{r}\sum_{j=2}^{n}\left( A_{1j}^{\alpha }\right) ^{2}.
\end{equation*}%
So we can write%
\begin{equation*}
Ric^{\ast }\left( X_{1}\right) \leq \frac{c+3}{2}\left( n-1\right) +\frac{3}{%
4}(c-1)\left\Vert CX_{1}\right\Vert ^{2}.
\end{equation*}%
Hence we obtain the following theorem:

\begin{theorem}
\label{thrm5} Let $\pi :M(c)\rightarrow B$ be an anti-invariant Riemannian
submersion from a Sasakian space form $\left( M(c),g\right) $ onto a
Riemannian manifold $\left( B,g^{^{\prime }}\right) $ such that $\xi $\ is
vertical. Then%
\begin{equation*}
Ric^{\ast }\left( X_{1}\right) \leq \frac{c+3}{4}\left( n-1\right) +\frac{3}{%
4}(c-1)\left\Vert CX_{1}\right\Vert ^{2}.
\end{equation*}%
The equality case of the inequality holds if and only if%
\begin{equation*}
A_{1j}=0,\text{ \ }j=2,...,n.
\end{equation*}
\end{theorem}

Since%
\begin{equation*}
2\tau =\underset{s=1}{\overset{n}{\sum }}Ric\left( X_{s},X_{s}\right) +%
\underset{k=1}{\overset{r}{\sum }}Ric\left( U_{k},X_{k}\right) ,
\end{equation*}%
\newpage
\begin{equation*}
2\tau =\underset{j,k=1}{\overset{r}{\sum }}R\left(
U_{j},U_{k},U_{k},U_{j}\right) +\underset{i=1}{\overset{n}{\sum }}\underset{%
k=1}{\overset{r}{\sum }}R\left( X_{i},U_{k},U_{k},X_{i}\right)
\end{equation*}%
\begin{equation}
+\underset{i,s=1}{\overset{n}{\sum }}R\left( X_{i},X_{s},X_{s},X_{i}\right) +%
\underset{s=1}{\overset{n}{\sum }}\underset{j=1}{\overset{r}{\sum }}R\left(
U_{j},X_{s},X_{s},U_{j}\right) ,  \label{S2}
\end{equation}%
where $\tau $ is the scalar curvature of $M(c)$. Since $M(c)$ is a Sasakian
space form, using (\ref{S2}) and (\ref{curv1}), we find%
\begin{equation}
2\tau =\frac{c+3}{4}\left( r\left( r-1\right) +n\left( n-1\right)
+2nr\right) +\frac{c-1}{4}\left( 4\left( r-1\right) +n+3tr\phi B\right) .
\label{S1}
\end{equation}%
On the other hand, from the Gauss-Codazzi type equations (\ref{gauss1}), (%
\ref{gauss2}) and (\ref{gauss3}), we have
\begin{equation*}
2\tau =2\widehat{\tau }+2\tau ^{\ast }+r^{2}\left\Vert H\right\Vert
^{2}+\sum_{k,j=1}^{r}g\left( T_{U_{k}}U_{j},T_{U_{k}}U_{j}\right)
\end{equation*}%
\begin{equation*}
+3\sum_{i,s=1}^{n}g\left( A_{X_{i}}X_{s},A_{X_{i}}X_{s}\right) -\underset{i=1%
}{\overset{n}{\sum }}\underset{k=1}{\overset{r}{\sum }}g\left( \left( \nabla
_{X_{i}}T\right) _{U_{k}}U_{k},X_{i}\right)
\end{equation*}%
\begin{equation*}
+\sum_{i=1}^{n}\sum_{k=1}^{r}\left\{ g\left(
T_{U_{k}}X_{i},T_{U_{k}}X_{i}\right) -g\left(
A_{X_{i}}U_{k},A_{X_{i}}U_{k}\right) \right\} -\underset{s=1}{\overset{n}{%
\sum }}\underset{j=1}{\overset{r}{\sum }}g\left( \left( \nabla
_{X_{s}}T\right) _{U_{j}}U_{j},X_{s}\right)
\end{equation*}%
\begin{equation}
+\sum_{s=1}^{n}\sum_{j=1}^{r}\left\{ g\left(
T_{U_{j}}X_{s},T_{U_{j}}X_{s}\right) -g\left(
A_{X_{s}}U_{j},A_{X_{s}}U_{j}\right) \right\} .  \label{S3}
\end{equation}%
Using (\ref{T1}) and (\ref{deltaN}), we get%
\begin{equation*}
2\tau =2\widehat{\tau }+2\tau ^{\ast }+\frac{1}{2}r^{2}\left\Vert
H\right\Vert ^{2}-\frac{1}{2}\sum_{s=1}^{n}\left[
T_{11}^{s}-T_{22}^{s}-...-T_{rr}^{s}\right] ^{2}
\end{equation*}%
\begin{equation*}
-2\sum_{s=1}^{n}\sum_{j=2}^{r}\left( T_{1j}^{s}\right)
^{2}+2\sum_{s=1}^{n}\sum_{2\leq j<k\leq r}^{r}\left[ T_{jj}^{s}T_{kk}^{s}-%
\left( T_{jk}^{s}\right) ^{2}\right] +6\sum_{\alpha
=1}^{r}\sum_{s=2}^{n}\left( A_{1s}^{\alpha }\right) ^{2}
\end{equation*}%
\begin{equation*}
+6\sum_{\alpha =1}^{r}\sum_{2\leq i<s\leq n}\left( A_{is}^{\alpha }\right)
^{2}+\sum_{i=1}^{n}\sum_{k=1}^{r}\left\{ g\left(
T_{U_{k}}X_{i},T_{U_{k}}X_{i}\right) -g\left(
A_{X_{i}}U_{k},A_{X_{i}}U_{k}\right) \right\}
\end{equation*}%
\begin{equation*}
-2\delta \left( N\right) +\sum_{s=1}^{n}\sum_{j=1}^{r}\left\{ g\left(
T_{U_{j}}X_{s},T_{U_{j}}X_{s}\right) -g\left(
A_{X_{s}}U_{j},A_{X_{s}}U_{j}\right) \right\} .
\end{equation*}%
By making use of (\ref{T1a}), (\ref{T1aa}) and (\ref{S1}) in the last
equality, we obtain%

\newpage

\begin{equation*}
\frac{c+3}{2}nr+\frac{c-1}{2}\left( 3\left( r-1\right) -n\right)
\end{equation*}%
\begin{equation*}
+2\underset{k=1}{\overset{r}{\sum }}R\left( U_{1},U_{k},U_{k},U_{1}\right) +2%
\underset{s=1}{\overset{n}{\sum }}R\left( X_{1},X_{s},X_{s},X_{1}\right)
\end{equation*}%
\begin{equation*}
=2\widehat{Ric}\left( U_{1}\right) +2Ric^{\ast }\left( X_{1}\right) +\frac{1%
}{2}r^{2}\left\Vert H\right\Vert ^{2}-\frac{1}{2}\sum_{s=1}^{n}\left[
T_{11}^{s}-T_{22}^{s}-...-T_{rr}^{s}\right] ^{2}
\end{equation*}%
\begin{equation*}
-2\sum_{s=1}^{n}\sum_{j=2}^{r}\left( T_{1j}^{s}\right) ^{2}+6\sum_{\alpha
=1}^{r}\sum_{s=2}^{n}\left( A_{1s}^{\alpha }\right)
^{2}+\sum_{i=1}^{n}\sum_{k=1}^{r}\left\{ g\left(
T_{U_{k}}X_{i},T_{U_{k}}X_{i}\right) -g\left(
A_{X_{i}}U_{k},A_{X_{i}}U_{k}\right) \right\}
\end{equation*}%
\begin{equation*}
-2\delta \left( N\right) +\sum_{s=1}^{n}\sum_{j=1}^{r}\left\{ g\left(
T_{U_{j}}X_{s},T_{U_{j}}X_{s}\right) -g\left(
A_{X_{s}}U_{j},A_{X_{s}}U_{j}\right) \right\} .
\end{equation*}%
We denote
\begin{equation*}
\left\Vert T^{V}\right\Vert ^{2}=\sum_{i=1}^{n}\sum_{k=1}^{r}g\left(
T_{U_{k}}X_{i},T_{U_{k}}X_{i}\right)
\end{equation*}%
and%
\begin{equation*}
\left\Vert A^{H}\right\Vert ^{2}=\sum_{i=1}^{n}\sum_{k=1}^{r}g\left(
A_{X_{i}}U_{k},A_{X_{i}}U_{k}\right) ,
\end{equation*}%
(see \cite{Gul-17}).

Since $\left( M(c),g\right) $ is a Sasakian space form, from (\ref{curv1}),
we obtain the following theorem:

\begin{theorem}
Let $\pi :M(c)\rightarrow B$ be an anti-invariant Riemannian submersion from
a Sasakian space form $\left( M(c),g\right) $ onto a Riemannian manifold $%
\left( B,g^{^{\prime }}\right) $ such that $\xi $\ is vertical. Then%
\begin{equation*}
\frac{c+3}{4}\left\{ nr+n+r-2\right\} +\frac{c-1}{4}\left\{ 3r-4-n\right.
\end{equation*}%
\begin{equation*}
\left. -\left( r-2\right) \eta \left( U_{1}\right) ^{2}+3\left\Vert
CX_{1}\right\Vert ^{2}\right\} \leq \widehat{Ric}\left( U_{1}\right)
+Ric^{\ast }\left( X_{1}\right) +\frac{1}{4}r^{2}\left\Vert H\right\Vert ^{2}
\end{equation*}%
\begin{equation*}
+3\sum_{\alpha =1}^{r}\sum_{s=2}^{n}\left( A_{1s}^{\alpha }\right)
^{2}-\delta \left( N\right) +\left\Vert T^{V}\right\Vert ^{2}-\left\Vert
A^{H}\right\Vert ^{2}.
\end{equation*}%
The equality case of the inequality holds if and only if%
\begin{equation*}
T_{11}^{s}=T_{22}^{s}+...+T_{rr}^{s},
\end{equation*}%
\begin{equation*}
T_{1j}=0,\text{ \ }j=2,...,r.
\end{equation*}
\end{theorem}

\textbf{Case II: }Assume that $\xi $\ is horizontal.

From (\ref{R1}), similar to Theorem \ref{thrm4}, we can state the following
theorem:

\begin{theorem}
\label{thrm44} Let $\pi :M(c)\rightarrow B$ be an anti-invariant Riemannian
submersion from a Sasakian space form $\left( M(c),g\right) $ onto a
Riemannian manifold $\left( B,g^{^{\prime }}\right) $ such that $\xi $\ is
horizontal. Then%
\begin{equation*}
\widehat{Ric}\left( U_{1}\right) \geq \frac{c+3}{4}\left( r-1\right) -\frac{1%
}{4}r^{2}\left\Vert H\right\Vert ^{2}.
\end{equation*}%
The equality case of the inequality holds if and only if%
\begin{equation*}
T_{11}^{s}=T_{22}^{s}+...+T_{rr}^{s},
\end{equation*}%
\begin{equation*}
T_{1j}=0,\text{ \ }j=2,...,r.
\end{equation*}
\end{theorem}

From (\ref{R2}), similar to Theorem \ref{thrm5}, we have the following
theorem:

\begin{theorem}
\label{thrm55} Let $\pi :M(c)\rightarrow B$ be an anti-invariant Riemannian
submersion from a Sasakian space form $\left( M(c),g\right) $ onto a
Riemannian manifold $\left( B,g^{^{\prime }}\right) $ such that $\xi $\ is
horizontal. Then%
\begin{equation*}
Ric^{\ast }\left( X_{1}\right) \leq \frac{c+3}{4}\left( n-1\right) +\frac{c-1%
}{4}\left\{ \left( 2-n\right) \eta \left( X_{1}\right) ^{2}-1+3\left\Vert
CX_{1}\right\Vert ^{2}\right\} .
\end{equation*}%
The equality case of the inequality holds if and only if%
\begin{equation*}
A_{1j}=0,\ j=2,...,n.
\end{equation*}
\end{theorem}

Since $\xi $\ is horizontal, from (\ref{S2}), we find%
\begin{equation*}
2\tau =\frac{c+3}{4}\left[ r\left( r-1\right) +n\left( n-1\right) +2nr\right]
+\frac{c-1}{4}\left[ n+3tr\phi B+4r-7\right] .
\end{equation*}%
Using the above equation, (\ref{S3}), (\ref{T1}), (\ref{T1a}), (\ref{T1aa})
and (\ref{deltaN}), we get%
\begin{equation*}
\frac{c+3}{2}nr+\frac{c-1}{2}\left( 2r-3\right)
\end{equation*}%
\begin{equation*}
+2\underset{k=1}{\overset{r}{\sum }}R\left( U_{1},U_{k},U_{k},U_{1}\right) +2%
\underset{s=1}{\overset{n}{\sum }}R\left( X_{1},X_{s},X_{s},X_{1}\right)
\end{equation*}%
\begin{equation*}
=2\widehat{Ric}\left( U_{1}\right) +2Ric^{\ast }\left( X_{1}\right) +\frac{1%
}{2}r^{2}\left\Vert H\right\Vert ^{2}-\frac{1}{2}\sum_{s=1}^{n}\left[
T_{11}^{s}-T_{22}^{s}-...-T_{rr}^{s}\right] ^{2}
\end{equation*}%
\begin{equation*}
-2\sum_{s=1}^{n}\sum_{j=2}^{r}\left( T_{1j}^{s}\right) ^{2}+6\sum_{\alpha
=1}^{r}\sum_{s=2}^{n}\left( A_{1s}^{\alpha }\right) ^{2}-2\delta \left(
N\right)
\end{equation*}%
\begin{equation*}
+\sum_{s=1}^{n}\sum_{j=1}^{r}\left\{ g\left(
T_{U_{j}}X_{s},T_{U_{j}}X_{s}\right) -g\left(
A_{X_{s}}U_{j},A_{X_{s}}U_{j}\right) \right\}
\end{equation*}%
\begin{equation*}
+\sum_{i=1}^{n}\sum_{k=1}^{r}\left\{ g\left(
T_{U_{k}}X_{i},T_{U_{k}}X_{i}\right) -g\left(
A_{X_{i}}U_{k},A_{X_{i}}U_{k}\right) \right\} .
\end{equation*}%
Hence in view of (\ref{curv1}), we obtain the following theorem:

\begin{theorem}
Let $\pi :M(c)\rightarrow B$ be an anti-invariant Riemannian submersion from
a Sasakian space form $\left( M(c),g\right) $ onto a Riemannian manifold $%
\left( B,g^{^{\prime }}\right) $ such that $\xi $\ is horizontal. Then%
\begin{equation*}
\frac{c+3}{4}\left\{ nr+n+r-2\right\} +\frac{c-1}{4}\left\{ 2r-4-\left(
n-2\right) \eta \left( X_{1}\right) ^{2}\right.
\end{equation*}%
\begin{equation*}
\left. +3\left\Vert CX_{1}\right\Vert ^{2}\right\} \leq \widehat{Ric}\left(
U_{1}\right) +Ric^{\ast }\left( X_{1}\right) +\frac{1}{4}r^{2}\left\Vert
H\right\Vert ^{2}
\end{equation*}%
\begin{equation*}
+3\sum_{\alpha =1}^{r}\sum_{s=2}^{n}\left( A_{1s}^{\alpha }\right)
^{2}-\delta \left( N\right) +\left\Vert T^{V}\right\Vert ^{2}-\left\Vert
A^{H}\right\Vert ^{2}.
\end{equation*}%
The equality case of the inequality holds if and only if%
\begin{equation*}
T_{11}^{s}=T_{22}^{s}+...+T_{rr}^{s},
\end{equation*}%
\begin{equation*}
T_{1j}=0,\text{ \ }j=2,...,r.
\end{equation*}
\end{theorem}

\bigskip

\end{document}